\documentclass{amsart}

\usepackage{color}
\usepackage{graphicx} 
\usepackage{latexsym} 
\usepackage{hyperref}
\usepackage{amsfonts} 
\usepackage[all]{xy} 
\usepackage{amssymb, mathrsfs, amsfonts, amsmath} 
\usepackage{amsbsy} 
\usepackage{amsfonts}
\usepackage{textcomp}
\usepackage{palatino}

\newtheorem*{theorem*}{Theorem}

\newtheorem{remark}{Remark}
 
\newtheorem{theorem}{Theorem} 
\newtheorem{example}{Example} 
\numberwithin{equation}{section}

\begin{document} 
	
\title[Recurrence and transience]{Recurrence and transience for normally reflected Brownian motion in warped product manifolds}

\author{Levi Lopes de Lima}
\address{Universidade Federal do Cear\'a,
		Departamento de Matem\'atica, Campus do Pici, Av. Humberto Monte, s/n, 60455-760,
		Fortaleza/CE, Brasil.}
\email{levi@mat.ufc.br}

\thanks{ The author was partially supported by CNPq/Brazil.}



\begin{abstract} 
We establish an integral test describing  the exact cut-off between recurrence and transience for  normally reflected Brownian motion in certain unbounded domains in a class of warped product manifolds. Besides extending a previous result by R. Pinsky \cite{P1}, who treated the case in which the ambient space is flat, our result recovers the classical test for the standard Brownian motion in model spaces. Moreover, it allows us to discuss the recurrence/transience  dichotomy for certain generalized tube domains around totally geodesic submanifolds in hyperbolic space. 
\end{abstract}

\maketitle

\section{Introduction and statement of the main result}\label{intro}

It is well-known that Brownian motion exhibits a sharp dichotomy regarding its long term behavior: either the typical path returns to any compact domain of the state space infinitely often (recurrence) or it escapes toward infinity, eventually leaving behind any bounded region (transience). In the geometric context  of  Brownian motion in a complete Riemannian manifold, several tests have been designed to single out which case actually occurs; see \cite{G} for a rather complete survey on this subject. It turns out that a  challenging problem here is  to determine, for a given family of metrics, the exact cut-off between the regimes. For instance, a classical result described in \cite[Proposition 3.1]{G} settles the problem for model spaces (rotationally symmetric manifolds). 

When the given manifold carries a non-empty boundary, it is natural to consider normally  reflected Brownian motion; we refer to \cite{SV} \cite{LS} \cite{IW} \cite{BH} \cite{H} \cite{P1} \cite{W} \cite{dL} for the construction and basic properties of this  diffusion process. In this setting, a question posed in \cite[Problem 6]{G} suggests an integral test to find out the precise threshold between recurrence and transience for a  generalized tube domain around a totally geodesic plane in flat space under the assumption that the  corresponding profile function does not oscillate too much at infinity. A rather satisfactory answer to this question was provided by R. Pinsky \cite{P1}, who based his argument on the construction of suitable Lyapunov functions. 
  
The aim of this note is to point out that the method in \cite{P1} can be adapted to handle certain unbounded domains in a class of warped product Riemannian manifolds displaying an $O_p\times O_q$ symmetry. In fact, from the argument presented in Section \ref{theproof} it transpires that the class of $O_p\times O_q$-invariant metrics considered here is the largest one to which the technique in \cite{P1} applies. Besides extending \cite[Theorem 1]{P1} to those metrics, our result retrieves the classical recurrence/transience test for the standard Brownian motion in model spaces referred to above. More importantly, as a further  application we will be  able to  discuss the recurrence/transience dichotomy for certain generalized tube domains around totally geodesic submanifolds in hyperbolic space whose profile functions decay exponentially (see Example \ref{slab} below). 

We now describe our main result (Theorem \ref{main} below). Fix an integer $q\geq 1$ and on the manifold  $M_q=(0,+\infty)\times \mathbb S^{q-1}$ consider the rotationally symmetric metric $k=dr^2+\xi^2(r)d\sigma_{q-1}^2$. Here, $d\sigma_{q-1}^2$ is the standard round metric on the unit sphere $\mathbb S^{q-1}\subset\mathbb R^q$ and $\xi:(0,+ \infty)\to\mathbb R$ is a smooth positive function. We assume that $k$ extends smoothly accross $r=0$ so as to define a Riemannian structure in $\mathbb R^q$. Also, we require that $\ddot\xi\geq 0$ everywhere, where the dot means derivative with respect to $r$. Note that this implies $\dot{\xi}\geq 1$. 
For $p\geq 2$ define on $M_p\times \mathbb R^q$ the 
metric $g=h+\phi^2(s)k$, where $h=ds^2+\theta^2(s)d\sigma_{p-1}^2$ 
is a complete rotationally symmetric metric in $\mathbb R^p$.
Here, $\theta,\phi:(0,+\infty)\to\mathbb R$ are smooth positive  functions.  Even though a version of our main result holds true for certain domains in warped product manifolds carrying an inner boundary (see Example \ref{schw} below), in this Introduction we take for granted that $g$ extends smoothly accross $s=0$, thus defining a Riemannian structure in $\mathbb R^{p+q}=\mathbb R^p\oplus\mathbb R^q$ which is obviously invariant under the standard $O_p\times O_q$-action. This invariance plays a key role in the computations leading to the proof of Theorem \ref{main}.

Now take a smooth positive function $f:[0,+\infty)\to \mathbb R$ and consider the domain $D_f=\{(x,y)\in \mathbb R^{p+q}; r<f(s)\}$, where $s=|x|$ and $r=|y|$. By construction, the Riemannian domain $(D_f,g)$ inherits the invariance under the $O_p\times O_q$-action mentioned above. We are interested in discussing the recurrence/transience dichotomy for normally reflected Brownian motion in this kind of domain. Our test involves the integral 
\[
I=\int^{\infty}\theta^{1-p}(s)\phi^{-q}(s)\xi(f)^{-q}(s)ds, 
\] 
which clearly depends only on the asymptotic geometry of $(D_f,g)$. Similarly to Pinsky's, our test requires suitable decay assumptions preventing  $f$ from oscillating too much at infinity. Here and in the following, primes denote ordinary or partial derivation with respect to $s$. 

\begin{theorem}
	\label{main}
	With the notation above assume that:
	\begin{itemize}
		\item $\phi^2(s)\xi(f)(s)f''(s)\to 0$, $\phi(s)f'(s)\to 0$,   $\phi'(s)\xi(f)(s)=O(1)$, $\dot{\xi}(f)(s)\to 1$, $\phi^2(s)\xi(f)(s)f'(s)=o(s)$ and $\xi(f)(s)\ddot{\xi}(f)(s)=O(1)$ as $s\to +\infty$;
		\item  $\int^{\infty}\phi^2(s)\frac{f'(s)^3}{\xi(f)(s)}ds$,  $\int^{\infty}\phi^2(s)\left(\frac{\theta'(s)}{\theta(s)}+
		\frac{\phi'(s)}{\phi(s)}\right)f'(s)^2ds$,
			$\int^{\infty}\xi(f)(s)\ddot{\xi}(f)(s)ds$,\,\,\,	  
		 $\int^{\infty}\left(\frac{f'(s)}{\xi(f)(s)}\left(1-(\dot{\xi}(f)(s))^2\right)
		\right) ds$ and 
		$\int^{\infty}\phi'(s)f'(s)ds$
		 are finite.  
	\end{itemize} 
	Then the following assertions hold:
	\begin{enumerate}
		\item If $I=+\infty$ and $ D\subset D_f$ then  normally reflected Brownian motion in $ D$ is recurrent;
		\item If $I<+\infty$ and $D_f\subset  D$ then  normally reflected Brownian motion in $ D$ is transient.
	\end{enumerate}
\end{theorem}

\begin{remark}
	\label{greg}{\rm As mentioned above, Theorem \ref{main} addresses the general question suggested by \cite[Problem 6]{G} and \cite[Theorem 1]{P1}. The result confirms that, at least for the class of warped domains considered here, there exists an integral test detecting the recurrence/trans\-ience cut-off if $f$ is ``slowly'' changing at infinity. However, as already pointed out in \cite{G},  an effective  test in case $f$ oscillates too much seems to be much harder to figure out.}
\end{remark}


\begin{remark}
	\label{phiconst}
{\rm
 The conclusions in the theorem can be rephrased in potential theoretic terms. Thus, recurrence means that any bounded subharmonic function $v$ on $D$ satisfying $\partial v/\partial \nu\geq 0$ along $\partial D$, where $\nu$ is the inward unit normal, is necessarily constant. By contrast, transience means that such non-constant functions exist; see \cite[Proposition 5.1]{G}.
}
\end{remark}

\begin{remark}
	\label{conser}{\rm In order to avoid explicitly mentioning explosion times, normally reflected Brownian motion $b_t$ in $(D_f,g)$ is tacitly assumed to be conservative in Theorem \ref{main}. Of course, recurrence already implies conservativeness. Moreover, in most examples of transience occurring in Section \ref{examp}, $(D_f,g)$ has the property that the Ricci tensor ${\rm Ric}$ and the second fundamental form $\mathcal A$ are both uniformly bounded from below. Thus, we may start with the Feyman-Kac formula for differential $1$-forms on manifolds with boundary \cite{H} \cite{dL} and proceed as in the proof of \cite[Theorem 3.2.6]{BGL} to make sure that conservativeness holds in those cases. For the sake of completeness we include here the well-known argument. If $\omega$ is a compactly supported $1$-form on $D_f$ satisfying absolute boundary conditions then the Feyman-Kac formula implies the estimate
	\begin{equation}\label{stocest2}
	|(e^{-\frac{1}{2}t\Delta^{(1)}_g}\omega)(b_0)|\leq \mathbb E_{b_0}\left(|\omega(b_t)|\exp\left(
	-\frac{1}{2}\int_0^t r(b_s)ds-\int_0^ta(b_s)dl_s
	\right)\right),\quad t>0.
	\end{equation}
	Here, $\Delta^{(1)}_g$ is the Hodge  Laplacian acting on $1$-forms, $r(x)=\inf_{|\omega|=1}{\langle\rm Ric}_x\omega,\omega\rangle$, $a(x)=\inf_{|\omega|=1}\langle \mathcal A_x\omega,\omega\rangle$
	and  
	$l_t$ is the boundary local time. 
	Now let $\varphi$ and $\eta$ be compactly supported functions on $N$ with $\partial \varphi/\partial \nu=0$ along $\partial N$ and $\eta=0$ in a neighborhood of $\partial N$. If $\Delta_g$ is the Laplacian acting on functions, a computation shows that  
	\[
	\int_{D_f}(e^{\frac{1}{2}t\Delta_g}\varphi-\varphi)\eta\, d{\rm vol}_g=-\frac{1}{2}\int_0^t\int_{D_f}\langle e^{-\frac{1}{2}\tau\Delta^{(1)}_g} d\varphi,d\eta\rangle\, d{\rm vol}_gd\tau, 
	\]
	and using (\ref{stocest2}) with $\omega=d\varphi$ we find that
	\[
	\left|	\int_{D_f}(e^{\frac{1}{2}t\Delta_g}\varphi-\varphi)\eta\, d{\rm vol}_g\right|\leq \frac{1}{2}|d\varphi|_\infty|d\eta|_1\sup_{0\leq \tau\leq t}e^{-{\underline r}\frac{\tau}{2}-\underline a\int_0^\tau dl^s},
	\]
	where $\underline r=\inf r$ and $\underline a=\inf a$  are both finite by assumption.
	Letting $\varphi$ approach $\bf 1$, the function identically equal to $1$, we conclude that $e^{\frac{1}{2}t\Delta_g}{\bf 1}={\bf 1}$ and conservativeness is established.	
		}
\end{remark}

\section{Some examples}\label{examp}

We now discuss a few interesting examples to which  Theorem \ref{main} applies.

\begin{example}
	\label{example2}{\rm If $\theta(s)=s$, $\xi(r)=r$ and $\phi=1$ we have
		\[
		I=\int^{\infty}s^{1-p}f^{-q}(s)ds,
		\]
		which yields \cite[Theorem 1]{P1}. 
		As discussed there, if applied to $f(s)=(1+s)^\alpha$, $\alpha<1$, this yields recurrence if $\alpha\leq\frac{2-p}{q}$ and transience if $\alpha>\frac{2-p}{q}$. Note that the volume of $D_f$ is finite if and only if $\alpha<-\frac{p}{q}$,
		in which case normally reflected Brownian motion is actually positive recurrent in the sense that it possesses an invariant probability measure \cite[Proposition 1]{P1}. 
		In particular, infinite volume coexists with null recurrence if and only if $-\frac{p}{q}\leq \alpha\leq \frac{2-p}{q}$. This should be compared to the results in Example \ref{slab} below, where the corresponding interval degenerates into a point.
		}
\end{example}

\begin{example}
	\label{example1} 
	{\rm If $\xi(r)=r$ and $\phi=f=1$ we have
	\[
	I=\int^{\infty}\theta^{1-p}(s)ds.
	\]
	Here, $(D_f,g)$ is the Riemannian product of the rotationally symmetric manifold $(\mathbb R^p,h)$ and the unit disk of dimension $q$ endowed with the flat metric. 
    After decomposing the normally reflected Brownian motion in $D_f$ as $b_t=(x_t,y_t)$, it is clear that $x_t$ is the  standard Brownian motion in  $(\mathbb R^p,h)$ (which we assume to be conservative; see Remark \ref{conser}). In this way   
	we recover the classical recurrence/tran\-sience test for model spaces \cite[Proposition 3.1]{G}.}
\end{example}

These examples illustrate the unifying character of Theorem \ref{main}. In addition, the test is flexible enough to allow for a discussion of the dichotomy in other interesting geometric backgrounds, as highlighted by the examples below. 

\begin{example}
	\label{slab}
	{\rm (Generalized tube domains in hyperbolic space)
				Take  
		$\theta(s)=\sinh s$, $\xi(r)=\sinh r$ and $\phi(s)=\cosh s$ so that   $h$ and $g$ are hyperbolic metrics in dimensions $p$ and $p+q$, respectively. Hence, $D_f$ is a generalized tube domain around a totally geodesic submanifold of dimension $p$.  
		Here, 
		\[
		I\sim\int^\infty e^{(1-p-q)s}\sinh^{-q}f(s)ds,
		\]
		for any $f$ as in Theorem \ref{main}.
		We apply the test to $f(s)=e^{\alpha s}$, $\alpha\in\mathbb R$. If $\alpha<0$ then  $f(s)\to 0$ uniformly as $s\to+\infty$, so that $\sinh f(s)=f(s)+o(e^{2\alpha s})$. Hence, when checking the assumptions in the theorem and figuring out the asymptotic shape of $I$ we may replace $\sinh f$ by $f$.   
		We easily see that the assumptions are satisfied if $\alpha<-1$ and that   
		\[
		I\sim
		\int^{\infty}e^{(1-p-q-\alpha q)s}ds.
		\] 
		Thus, recurrence takes place if $\alpha\leq \frac{1-p-q}{q}$ while transience occurs if $\alpha>\frac{1-p-q}{q}$. 
		We note that our test does not directly apply if $\alpha\geq -1$, since in this case $f$ oscillates ``too much'' at infinity and the assumptions in the theorem fail to hold; see Remark \ref{greg}. However, since  $\frac{1-p-q}{q}<-1$, we can appeal to domain monotonicity 
		(see the discussion following \cite[Proposition 1]{P1})
		to infer that transience also holds in this range of $\alpha$. 
		Notice that if $\alpha\geq 0$ this also follows from the fact that  $D_f$ contains the slab corresponding to the function identically equal to $1$, where transience should hold by a variant of the argument in Example \ref{example1}.
		We also observe that  $D_f$ has finite volume if and only if  $\alpha<\frac{1-p-q}{q}$. Thus, infinite volume coexists with null recurrence if and only if $\alpha=\frac{1-p-q}{q}$; compare with Example \ref{example2}. We remark that  if $f$ has superlinear exponential decay, $f(s)=e^{\alpha s^\gamma}$, $\alpha<0$, $\gamma>1$, then one always has recurrence. This can be seen  either by domain monotonicity (just compare with the case of liner exponential decay analysed above) or directly from the test (the assumptions in Theorem \ref{main} are always satisfied and one finds that $I\sim\int^{\infty}e^{(1-p-q)s-\alpha\gamma e^{\alpha s^\gamma}}ds=+\infty$). Finally, domain monotonicity also leads to transience in the sublinear case $\gamma<1$.} 
		\end{example}
	
\begin{example}\label{schw}{\rm (Generalized tube domains in Schwarzschild-type space)
Take parameters $m>0$, $\epsilon=0,-1$, and let $\mathfrak r_{m,\epsilon}$ be the unique positive zero of 
\[
\zeta_{m,\epsilon}(\mathfrak r)={1-\epsilon \mathfrak r^2-\frac{2m}{\mathfrak r^{p-2}}}.
\]
On $(\mathfrak r_{m,\epsilon},+\infty)\times \mathbb S^{p-1}\times \mathbb R^q$ define the metric $g_{m,\epsilon}=h_{m,\epsilon}+\zeta_{m,\epsilon}(\mathfrak r)k_\epsilon$, 
where $k_\epsilon=dr^2+\vartheta_\epsilon^2(r)\sigma_{q-1}^2$, with $\vartheta_0(r)=r$, $\vartheta_{-1}(r)=\sinh r$, and 
\[
h_{m,\epsilon}=\frac{d\mathfrak r^2}{\zeta_{m,\epsilon}(\mathfrak r)}+\mathfrak r^2d\sigma^2_{p-1}. 
\]
Even though this metric has many interesting geometric properties (for instance, if $q=1$ a well-known calculation shows that it is Einstein, ${\rm Ric}_{g_{m,\epsilon}}=p\epsilon g_{m,\epsilon}$), it is far from having constant sectional curvature. To probe its structure near the inner boundary $\mathfrak r=\mathfrak r_{m,\epsilon}$, introduce 
the radial parameter $s=\int_{\mathfrak r_{m,\epsilon}} {d\mathfrak r}/{\sqrt{\zeta_{m,\epsilon}(\mathfrak r)}}$. It follows that 
$h_{m,\epsilon}=ds^2+\theta^2(s)d\sigma_{p-1}^2$, where $\theta(0)=\theta'(0)=0$, which proves that $h_{m,\epsilon}$ extends smoothly to $[\mathfrak r_{m,\epsilon},+\infty)\times \mathbb S^{p-1}$ so that $\mathfrak r=\mathfrak r_{m,\epsilon}$ is a totally geodesic sphere. 
On the other hand, $g_{m,\epsilon}=h_{m,\epsilon}+\phi^2(s)k_\epsilon$, where ${\phi}(0)=0$, which means that the inner boundary (horizon) is null in the sense that the metric degenerates in the vertical directions as $s\to 0$ (equivalently, $\mathfrak r\to\mathfrak r_{m,\epsilon}$). However, if we choose $\underline s>0$ and consider $f:[\underline s,+\infty)\to[0,+\infty)$  with $f(\underline s)=0$ and $f(s)>0$ if $s>\underline s$ then $\overline D_f$ never intersects the null boundary and the argument leading to Theorem \ref{main} carries over to this situation. Since $g_{m,\epsilon}$ approaches the flat or hyperbolic metric as $s\to+\infty$ (according to  $\epsilon$ being $0$ or $-1$)   
and our test only involves the asymptotic geometry of $(D_f,g_{m,\epsilon})$, the results in Examples \ref{example2} and \ref{slab} above hold for the corresponding generalized tube domains in this Schwarzschild-type space.} 
\end{example}

\section{The proof of Theorem \ref{main}}\label{theproof}

The proof of Theorem \ref{main} follows the reasoning in \cite{P1} closely and relies on a well-known test for recurrence/transience based on the construction of Lyapunov functions \cite[Chapter 6]{P2}. Here we merely indicate the required modifications to carry over the argument to our more general setting. These are mainly due to the presence of the warping functions $\phi$ and $\xi$. We note that by domain monotonicity it suffices to prove the result for $D=D_f$.  

We retain the notation of the Introduction and consider the domain $D_f\subset\mathbb R^{p+q}$ defined by the condition $r<f(s)$ and endowed with the warped metric $g=h+\phi^2k$. If ${\rm e}_y=y/|y|$ then the inward unit normal vector on $\partial D_f$ is 
 \[
\nu=\frac{1}{W}\left({\phi f'\partial_s-\phi^{-1}{\rm e}_y}\right),\quad W= \sqrt{1+\phi^2f'^2}.
\]
The first step in the argument is an observation which follows readly from It\^o's formula: if $b_t$ is normally reflected Brownian motion in $D_f$ then for any smooth function $u:\overline D_f\to\mathbb R$ with $\langle \nabla u,\nu\rangle=0$ along $\partial \overline D_f$ the process
\begin{equation}\label{ito}
\Pi^u_t=u(b_t)-\frac{1}{2}\int_0^t\Delta_gu(b_\tau)d\tau
\end{equation}
is a local martingale. 
To exploit this we follow \cite{P1} and seek Lyapunov functions 
$u^\pm=u^\pm(\rho,r)$  as follows. We consider $(s,r)$ such that $0\leq r\leq f(s)$ and set
\[
\rho=F(s,r)=-\frac{1}{2}\phi^2(s)\frac{\xi(f)'(s)}{\xi(f)(s)\dot{\xi}(f)^2(s)}\xi(r)^2+s+\frac{1}{2}\phi^2(s)\frac{\xi(f)'(s)\xi(f)(s)}{\dot{\xi}(f)^2(s)}.
\]
Thus, $u^\pm$ will be determined by the requirement that 
\[
u^\pm(\Psi(s,r))=\psi^\pm(s),
\]
where $\Psi(s,r)=(F(s,r),r)$ and $\psi^\pm:[0,+\infty)\to \mathbb R$ is a positive function to be chosen later. We observe that for each $s$ the map $r\mapsto \Psi(s,r)$ parametrizes a curve which happens to  be a level set of $u^\pm$. Since 
$\Psi(s,f(s))=(s,f(s))$, this curve always hits $\partial D_f$. Moreover, it does it orthogonally because  
\[
\partial_r\Psi(s,f(s))=-\phi^2f'\partial_s+{\rm e}_y=-\phi W\nu.
\] 
In particular, $\langle \nabla_gu^\pm,\nu\rangle=0$ along any $\partial \overline D_f\cap\{\rho\geq \rho_0\}$. Thus, by stopping the process $\Pi_t^{u^\pm}$ with the exit time of domains of the type $\{\rho\leq\rho_0\}$ and taking expectation we see from (\ref{ito}) that suitable Lyapunov functions can be constructed if we make sure that $\Delta_gu^\pm$ has a definite sign, possibly vanishing somewhere, and control the growth of $u^\pm$ as $\rho\to+\infty$.

In the following calculations we sometimes suppress the dependence of the various functions on their arguments since this will cause no confusion. Accordingly, we define 
\begin{equation}\label{defL}
L=\frac{\xi(f)'}{\xi(f)\dot{\xi}(f)^2}=\frac{f'}{\xi(f)\dot{\xi}(f)},
\end{equation}
\[
E=(\phi^2L)', \quad N=\frac{\phi^2f'\xi(f)}{\dot{\xi}(f)}, 
\]
and 
\[
C=\rho'=1-\frac{1}{2}E\xi(f)^2+\frac{1}{2}N'.
\]
By setting $u=u^\pm$, $\psi=\psi^\pm$ for simplicity and taking partial derivatives of $u(\rho,r)=\psi(s)$ up to second order we obtain
\[
-\phi^2L\xi\dot{\xi} u_\rho+u_r=0,
\]
\[
\phi^4L^2\xi^2\dot{\xi}^2u_{\rho\rho}-\phi^2L(\dot{\xi}^2+\xi\ddot{\xi}) u_\rho+u_{rr}-2\phi^2\xi\dot{\xi}Lu_{\rho r}=0, 
\]
\[
Cu_\rho=\psi',
\]
\[
C^2u_{\rho r}+ C'u_\rho=\psi'',
\]
and 
\[
Cu_{\rho r}-\phi^2\xi\dot{\xi} LCu_{\rho\rho}-E\xi\dot{\xi}u_\rho=0. 
\]
Since 
\begin{equation}\label{ellorbit}
\Delta_gu=u_{\rho\rho}+\left((p-1)\frac{ \theta_\rho}{\theta}+q\frac{\phi_{\rho}}{\phi}\right)u_\rho+\frac{1}{\phi^2}\left(u_{rr}+(q-1)\frac{\dot{\xi}}{\xi}u_r\right),
\end{equation}
a direct computation gives
\begin{equation}\label{laplace}
\Delta_gu = A\psi''+B\psi', 
\end{equation}
where 
\[
A=\frac{G}{C^2},\quad G=1+\phi^2L^2\xi^2\dot{\xi}^2,
\]
and 
\[
B=q\frac{L'}{C}\dot\xi^2+(p-1)\frac{1}{C}\frac{\theta_\rho}{\theta}+q\frac{1}{C}
\frac{\phi_{\rho}}{\phi}+\frac{\xi\ddot{\xi}}{C}-\frac{C'G}{C^3}+
2\frac{LE\xi^2\dot{\xi}^2}{C^2}. 
\]

The {\em rationale} behind the computation above can be described as follows. The invariance of $(D_f,g)$ under the $O_p\times O_q$-action on $\mathbb R^{p+q}$ suggests that the (reflected) Markov diffusion process on $D_f/O_p\times O_q=\{(s,r)\in\mathbb R^2; s\geq 0, r\geq 0, r<f(s)\}$ induced by the elliptic operator in the right-hand side of (\ref{ellorbit}) should contain a great deal of information regarding the original reflected Brownian motion in $(D_f,g)$. However, a direct analysis of this bidimensional diffusion is  out of the question. Fortunately, consideration of the prospective Lyapunov function $u=u^{\pm}$ allows us to perform a further dimensional reduction: 
for each $r$ the right-hand side of  (\ref{laplace})
defines a second order differential operator acting on the auxiliary function $\psi=\psi(s)$. Under our assumptions, for $s$ large enough this operator is elliptic with sufficiently regular coefficients. It is  well-known that the  recurrence/transience cut-off for the one-dimensional Markov diffusion process  generated by  such an operator involves the consideration of integrals of the type
\begin{equation}\label{inttype}
\int e^{-\int^s\frac{B}{A}(\rho)d\rho}ds;
\end{equation}
see \cite[Chapter VI]{IW} \cite[Chapter 5]{P2}. Thus, one is naturally led to 
squeeze  $B/A$ between functions depending only on $s$, and this is precisely where the decay conditions on $f$ come into play. We now explain how this can be done.

Indeed, after algebraic rearrangements we find that  
\begin{eqnarray*}
\frac{B}{A} & = & -\frac{C'}{C}+
\left(2-\frac{q}{2}\right)\frac{G'}{G} + (p-1)\frac{\theta'}{\theta}+q\frac{\phi'}{\phi}+q L +\frac{1}{2}q LN' \\
& & \quad 
+\left(\frac{q}{2}-2\right)\frac{\xi^2\phi^2 LL'}{H}
-(p-1)\frac{\theta'}{\theta}\frac{\xi^2\phi^2L^2}{H}\\
& & \quad\quad 
-q\frac{\phi'}{\phi}\frac{\xi^2\phi^2L^2}{H}
-q\frac{\xi^2\phi^2L^3
	\left(1+\frac{N'}{2}\right)}{H}+ \frac{C\xi\ddot{\xi}}{H\dot{\xi}^2}\\
& & \quad \quad \quad +qL\left(C+\frac{E\xi^2}{2}\right)\frac{1-\dot{\xi}^2}{H},
\end{eqnarray*}
where $H=\dot{\xi}^{-2}G$, $\theta'=\rho'\theta_\rho=C\theta_\rho$ and  similarly for $\phi'$.
For a fixed $s$ and as a function of $r\in [0,f(s)]$, the first four logarithm derivatives in the right-hand side above take their maximum and minimum values at the extremities of the interval. Thus, the functions
\[
\mathcal G^+(s)=\sup_{0\leq s\leq f(r)}\frac{G'}{G}(s,r), \quad \mathcal G^-(s)=\inf_{0\leq s\leq f(r)}\frac{G'}{G}(s,r),
\]
\[
\mathcal C^+(s)=\sup_{0\leq s\leq f(r)}\frac{C'}{C}(s,r), \quad \mathcal C^-(s)=\inf_{0\leq s\leq f(r)}\frac{C'}{C}(s,r),
\]
\[
\Theta^+(s)=\sup_{0\leq s\leq f(r)}\frac{\theta'}{\theta}(s,r), \quad  \Theta^-(s)=\inf_{0\leq s\leq f(r)}\frac{\theta'}{\theta}(s,r),
\]
and 
\[
\Phi^+(s)=\sup_{0\leq s\leq f(r)}\frac{\phi'}{\phi}(s,r), \quad  \Phi^-(s)=\inf_{0\leq s\leq f(r)}\frac{\phi'}{\phi}(s,r),
\]
satisfy 
\[
\int_{s_0}^t\mathcal G^\pm(\rho)d\rho=\log\frac{{G}(t,f_{G^\pm}(t))}{G(s_0,f_{G^\pm}(s_0)))}, \quad f_{G^\pm}(s)\in[0,f(s)], 
\]
\[
\int_{s_0}^t\mathcal C^\pm(\rho)d\rho=\log\frac{{C}(t,f_{C^\pm}(t))}{C(s_0,f_{C^\pm}(s_0)))}, \quad f_{C^\pm}(s)\in[0,f(s)], 
\]
\[
\int_{s_0}^t \Theta^\pm(\rho)d\rho=\log\frac{{\theta}(t,f_{\theta^\pm}(t))}{\theta(s_0,f_{\theta^\pm}(s_0)))}, \quad f_{\theta^\pm}(s)\in[0,f(s)], 
\]
and 
\[
\int_{s_0}^t \Phi^\pm(\rho)d\rho=\log\frac{{\phi}(t,f_{\phi^\pm}(t))}{\phi(s_0,f_{\phi^\pm}(s_0)))}, \quad f_{\phi^\pm}(s)\in[0,f(s)], 
\]

On the other hand, from the pointwise assumptions there exist positive constants $K_0,K_1,K_2, K_3$ such that
\[
\left|\frac{C\xi\ddot{\xi}}{H\dot{\xi}^2}+qL\left(C+\frac{E\xi^2}{2}\right)\frac{1-\dot{\xi}^2}{H}\right|
\leq K_0\left(\frac{|f'|}{\xi(f)}\left|1-\dot{\xi}(f)^2\right|+\xi(f)\ddot{\xi}(f)
\right),
\]
\[
\left(\frac{q}{2}-2\right)\left|\frac{\xi^2\phi^2LL'}{H}\right|\leq K_1\left(
\phi^2\frac{|f'|^3}{\xi(f)}+\phi^2|f'||f''|\right),
\]
\[
q\left|\frac{\xi^2\phi^2L'^3
	\left(1+\frac{N'}{2}\right)}{H}\right|\leq K_2\phi^2\frac{|f'|^3}{\xi(f)},
\]
and 
\[
(p-1)\left|\frac{\theta'}{\theta}\frac{\xi^2\phi^2L^2}{H}\right|
+
q\left|\frac{\phi'}{\phi}\frac{\xi^2\phi^2L^2}{H}\right|
\leq K_3\phi^2\left(\frac{|\theta'|}{\theta}+\frac{|\phi'|}{\phi}\right){f'^2}.
\]
Putting together these estimates, it follows from the expression for $B/A$ above that the functions
\begin{eqnarray}
\label{gammaf}
\Gamma^\pm & = & (p-1)\frac{\theta'}{\theta}+q\frac{\phi'}{\phi}+q \frac{f'}{\xi(f)\dot{\xi}(f)} + \frac{1}{2}q LN'\pm K_1
\phi^2|f'||f''|\nonumber \\
& & \quad \pm (K_1+K_2)\phi^2\frac{|f'|^3}{\xi(f)}\pm K_3\phi^2\left(\frac{|\theta'|}{\theta}+\frac{|\phi'|}{\phi}\right){f'^2}-\mathcal C^\pm+\left(2-\frac{q}{2}\right)\mathcal G^\pm\nonumber\\
& & \quad \quad \pm K_0\left(\frac{|f'|}{\xi(f)}\left|1-\dot{\xi}(f)^2\right|+\xi(f)\ddot{\xi}(f)
\right)
\end{eqnarray}
satisfy
\[
\Gamma^-(s)\leq \inf_{r\leq f(s)}\frac{B}{A}(r,s)\leq \sup_{r\leq f(s)}\frac{B}{A}(r,s)\leq \Gamma^+(s).
\]
As promised above, we have been able to control $B/A$ in terms of  quite explicit functions depending only on $s$. In light of  (\ref{inttype})  
we choose
\begin{equation}\label{expint}
\psi^\pm(s)=\int_{s_0}^se^{-\int_{s_0}^t\Gamma^{\pm}(\rho)d\rho}dt,\quad s_0>0,  
\end{equation}
which immediately gives $\Delta_gu^+ \leq 0$ and $\Delta_gu^-\geq 0$. Thus, the proof of Theorem \ref{main} will be completed as soon as we check that $\psi^+$ blows up at infinity whereas $\psi^-$ remains bounded. In the first case recurrence occurs while transience takes place in the second case; see \cite[Chapter 6]{P2}.

We claim that the contributions of the last six terms in the right-hand side of (\ref{gammaf}) to the exponential integrand in (\ref{expint}) are bounded and bounded away from $0$ for large $s_0$. This step of the proof uses the estimates above and the assumptions in Theorem 1 just as  in \cite{P1}, the only novelty coming from the term involving $f''$. But note that
\[
\int_{s_0}^t\phi^2f'f''ds=\frac{1}{2}(\phi^2f'^2)|_{s_0}^t-\int_{s_0}^t\phi\phi'f'^ 2ds, 
\]
and the claim follows since $\phi f'\to 0$. 
Also,   
\begin{eqnarray*}
\int_{s_0}^tLN'\,ds & = & 2\int_{s_0}^t\phi\phi'f'^2ds+\int_{s_0}^t\frac{\phi^2 f'f''}{\dot{\xi}(f)^2}ds\\
& & \quad 
+\int_{s_0}^t\phi^2\frac{f'^3}{\xi(f)\dot{\xi}(f)}ds-2\int_{s_0}^t
\phi^2\frac{f'^3\xi(f)\ddot{\xi}(f)}{\xi(f)\dot{\xi}(f)^2}ds. 
\end{eqnarray*}
Note that the contribution coming from the last two integrals in the right-hand side above remains likewise controlled. Moreover, taken together the first two integrals asymptote the quantity
\[
\phi^2f'^2|_{s_0}^t
-\int_{s_0}^t\phi^2f'f''ds,
\] 
so that the overall contribution coming from the term $qLN'/2$ is controlled as well.
Combining this with the estimates 
\[
(1-k_1(s_0))\frac{\theta(t)}{\theta(s_0)}\leq 
\frac{{\theta}(t,f_{\theta^\pm}(t))}{\theta(s_0,f_{\theta^\pm}(s_0)))}
\leq
(1+k_1(s_0))\frac{\theta(t)}{\theta(s_0)}
\]
and
\[
(1-k_2(s_0))\frac{\phi(t)}{\phi(s_0)}\leq 
\frac{{\phi}(t,f_{\theta^\pm}(t))}{\phi(s_0,f_{\phi^\pm}(s_0)))}
\leq
(1+k_2(s_0))\frac{\phi(t)}{\phi(s_0)},
\]
where $k_1(s_0),k_2(s_0)=o(s_0)$ as $s_0\to +\infty$,
we see that, up to a multiplicative factor that goes to $1$ as $s_0\to+\infty$,  
\[
e^{-\int_{s_0}^t\Gamma^\pm(\rho)d\rho}
\sim
\left(\frac{\theta(t)}{\theta(s_0)}\right)^{1-p}
\left(\frac{\phi(t)}{\phi(s_0)}\right)^{-q}
e^{-q\int_{s_0}^t L(\rho)d\rho}.
\]
From (\ref{defL}) and the assumption that $\dot{\xi}(f)(s)\to 1$ as $s\to +\infty$ it follows that 
\[
e^{-q\int_{s_0}^t L(\rho)d\rho}\sim
\left(\frac{\xi(f)(t)}{\xi(f)(s_0)}\right)^{-q},
\]
so that from (\ref{expint}) we conclude  that
$\psi^+$ blows up if $I=+\infty$ and $\psi^-$ remains bounded if $I<+ \infty$. As already observed, this completes  the proof of Theorem \ref{main}.


\begin{thebibliography}{999999}
	
\markboth{}{}


\bibitem[BGL]{BGL} Bakry, D., Gentil, I., Ledoux, M., {\em Analysis and geometry of Markov diffusion operators}. Grundlehren der Mathematischen Wissenschaften , 348. Springer, 2014.

\bibitem[BH]{BH} Bass, R. F., Hsu, P., The semimartingale structure of reflecting Brownian motion. {\em Proc. Amer. Math. Soc.} 108 (1990), no. 4, 1007-1010.

\bibitem[dL]{dL} de Lima, L. L., A Feynman-Kac formula for differential forms on manifolds with boundary and applications, {\em 	arXiv:1512.01153}.

\bibitem[G]{G}  Grigor{'}yan, A., Analytic and geometric background of recurrence and non-explosion of the Brownian motion on Riemannian manifolds. {\em Bull. Amer. Math. Soc. (N.S.)} 36 (1999), no. 2, 135-249.

\bibitem[H]{H} Hsu, E. P., Multiplicative functional for the heat equation on manifolds with boundary. {\em Michigan Math. J.} 50 (2002), no. 2, 351-367.

\bibitem[IW]{IW} Ikeda, N., Watanabe, S., {\em Stochastic differential equations and diffusion processes.}  North-Holland Publishing Co., Amsterdam; Kodansha, Ltd., Tokyo, 1989.


\bibitem[LS]{LS} Lions, P.-L., Sznitman, A.-S.,
Stochastic differential equations with reflecting boundary conditions.
{\em Comm. Pure Appl. Math.} 37 (1984), no. 4, 511-537.

\bibitem[P1]{P1} Pinsky, R. G.,
Transcience/recurrence for normally reflected Brownian motion in unbounded domains. 
{\em Ann. Probab.} 37 (2009), no. 2, 676-686.  

\bibitem[P2]{P2}  Pinsky, R. G., {\em Positive harmonic functions and diffusion}. Cambridge Studies in Advanced Mathematics, 45. Cambridge University Press, Cambridge, 1995. 

\bibitem[SV]{SV} Stroock, D. W., Varadhan, S. R. S.,
Diffusion processes with boundary conditions. 
{\em Comm. Pure Appl. Math.} 24 (1971) 147-225.

\bibitem[W]{W} Wang, F.-Y. {\em 
Analysis for diffusion processes on Riemannian manifolds.}
World Scientific Publishing Co. Pte. Ltd., Hackensack, NJ, 2014.

\end{thebibliography}
\end{document}